\documentclass[a4paper,11pt]{amsart}
\usepackage[OT2,OT1]{fontenc}
\usepackage[latin1]{inputenc}
\usepackage{amsmath}
\usepackage{amsthm}
\usepackage{amssymb}

\makeatletter
\renewcommand{\@evenhead}{\textsc{\small \thepage \hfil Peter Jossen and Antonella Perucca \hfil}}
\renewcommand{\@oddhead}{\textsc{\small \hfil A Counterexample to the local-global principle of detecting\ldots \hfil \thepage }}
\makeatother

\newtheorem*{cla}{Claim}

\DeclareMathOperator{\Mat}{Mat}
\DeclareMathOperator{\tr}{tr}

\newcommand{\fp}{\mathfrak p}

\newcommand{\IZ}{\mathbb Z}

\author{Peter Jossen and Antonella Perucca}
\title{A counterexample to the local-global principle of linear dependence for abelian varieties}
\date{}

\begin{document}

\begin{abstract}
Let $A$ be an abelian variety defined over a number field $k$. Let $P$ be a point in $A(k)$ and let $X$ be a subgroup of $A(k)$. 
Gajda in 2002 asked whether it is true that the point $P$ belongs to $X$ if and only if the point $(P \bmod \fp)$ belongs to $(X \bmod \fp)$ for all but finitely many primes $\fp$ of $k$. We answer negatively to Gajda's question. \end{abstract}

\maketitle

\vspace{4mm}
\begin{par}
Let $A$ be an abelian variety defined over a number field $k$. Let $P$ be a point in $A(k)$ and let $X$ be a subgroup of $A(k)$. Suppose that for all but finitely many primes $\fp$ of $k$ the point $(P \bmod \fp)$ belongs to $(X \bmod \fp)$. Is it true that $P$ belongs to $X$? This question, which was formulated by Gajda in 2002, was named the problem of detecting linear dependence. The problem was addressed in several papers (see the bibliography) but the question was still open. In this note we show that the answer to Gajda's question is negative by providing a counterexample.
\end{par}

\vspace{4mm}
\begin{par}
Let $k$ be a number field and let $E$ be an elliptic curve without complex multiplication over $k$ such that there are points $P_1$, $P_2$, $P_3$ in $E(k)$ which are $\IZ$--linearly independent. Let $P \in E^3(k)$ and $X \subseteq E^3(k)$ be the following:
$$P := \begin{pmatrix}  P_1\\ P_2\\ P_3 \end{pmatrix} \qquad\qquad X := \Big\{ MP \in E^3(k) \:\Big|\:\: M\in\Mat(3,\IZ), \tr M=0 \Big\}$$
So the group $X$ consists of the images of the point $P$ via the endomorphisms of $E^3$ with trace zero. Since the points $P_i$ are $\IZ$--independent and since $E$ has no complex multiplication, the point $P$ does not belong to $X$. Notice that no non-zero multiple of $P$ belongs to $X$.
\end{par}

\vspace{4mm}
\begin{cla}
Let $\fp$ be a prime of $k$ where $E$ has good reduction. The image of $P$ under the reduction map modulo $\fp$ belongs to the image of $X$.
\end{cla}

\vspace{4mm}
\begin{par}
For the rest of this note, we fix a prime $\fp$ of good reduction for $E$. We write $\kappa$ for the residue field of $k$ at $\fp$. To ease notation, we now let $E$ denote the reduction of the given elliptic curve modulo $\fp$ and write $P_1$, $P_2$, $P_3$, $P$  and $X$ for the image of the given points and the given subgroup under the reduction map modulo $\fp$.\\
Our aim is to find an integer matrix $M$ of trace zero such that $P = MP$ in $E^3(\kappa)$. 
\end{par}

\vspace{4mm}
\begin{par}
Let $\alpha_1$ be the smallest positive integer such that $\alpha_1 P_1$ is a linear combination of $P_2$ and $P_3$. Similarly define $\alpha_2$ and $\alpha_3$ for $P_2$ and $P_3$ respectively. There are integers $m_{ij}$ such that
\begin{eqnarray*}
\alpha_1 P_1 + m_{12}P_2+ m_{13}P_3 & = & O \\
m_{21}P_1+ \alpha_2 P_2 + m_{23}P_3 & = & O \\
m_{31}P_1 + m_{32}P_2+ \alpha_3 P_3 & = & O
\end{eqnarray*}
Assume that the greatest common divisor of $\alpha_1$, $\alpha_2$ and $\alpha_3$ is $1$ (we prove this assumption later). We can thus find integers $a_1, a_2, a_3$ such that 
$$3= \alpha_1 a_1 + \alpha_2 a_2 + \alpha_3a_3$$
Write $m_{ii} := 1-\alpha_ia_i$, so that in particular $m_{11}+m_{22}+m_{33} = 0$. Then we have
\begin{equation*}
\begin{pmatrix}
 m_{11} & -a_1 m_{12} & -a_1 m_{13}\\
 -a_2 m_{21} & m_{22} & -a_2 m_{23}\\
 -a_3 m_{31} & -a_3 m_{32} & m_{33} 
\end{pmatrix}
\begin{pmatrix} P_1\\ P_2\\P_3\end{pmatrix}=\begin{pmatrix}  P_1\\ P_2\\ P_3   \end{pmatrix}
\end{equation*} 
Notice that the above matrix has integer entries and trace zero. Hence we are left to prove that the greatest common divisor of $\alpha_1$, $\alpha_2$ and $\alpha_3$ is indeed $1$.
\end{par}

\vspace{4mm}
\begin{par}
Fix a prime number $\ell$ and let us show that $\ell$ does not divide $\gcd(\alpha_1,\alpha_2,\alpha_3)$. Suppose on the contrary that $\ell$ divides $\gcd(\alpha_1,\alpha_2,\alpha_3)$. This is equivalent to saying that $\ell$ divides all the coefficients appearing in any linear relation between $P_1$, $P_2$ and $P_3$. In particular, this implies that $\ell$ divides the order of $P_1$, $P_2$ and $P_3$ in $E(\kappa)$.
\end{par}

\vspace{4mm}
\begin{par}
It is well--known that the group $E(\kappa)[\ell]$ is either trivial, isomorphic to $\IZ/\ell\IZ$ or isomorphic to $(\IZ/\ell\IZ)^2$. In any case, the intersection $X \cap E(\kappa)[\ell]$ is generated by two elements or less. Without loss of generality, let us suppose that the subgroup of $X$ generated by $P_2$ and $P_3$ contains $X \cap E(\kappa)[\ell]$.\\
We are supposing that $\ell$ divides all the coefficients appearing in any linear relation of the points. Let $\alpha_1=x_1\ell$ and write $x_1\ell P_1+ x_2\ell P_2+ x_3\ell P_3= O$ for some integers $x_2$ and $x_3$. It follows that 
$$x_1P_1+x_2P_2+x_3P_3 =T$$
for some point $T$ in $X\cap E(\kappa)[\ell]$. The point $T$ is a linear combination of $P_2$ and $P_3$. This contradicts the minimality of $\alpha_1$.
\end{par}

\bibliographystyle{amsplain}
\nocite{*}

\providecommand{\bysame}{\leavevmode\hbox to3em{\hrulefill}\thinspace}
\providecommand{\MR}{\relax\ifhmode\unskip\space\fi MR }
\providecommand{\MRhref}[2]{%
  \href{http://www.ams.org/mathscinet-getitem?mr=#1}{#2}
}
\providecommand{\href}[2]{#2}

\vspace{15mm}
$$ $$
\hspace{20mm}
\begin{minipage}[]{120mm}
Peter Jossen\\
Department of Mathematics, Central European University\\
N\'ador u. 9, 1051 Budapest, Hungary\\[1mm]
{\tt peter.jossen@gmail.com}\\[8mm]
Antonella Perucca\\
Section des Math\'ematiques, \'Ecole Polytechnique F\'ed\'erale de Lausanne\\ 
EPFL, MA C3 604, Station 8, 1015 Lausanne, Switzerland\\[1mm]
{\tt antonella.perucca@epfl.ch}
\end{minipage}

\end{document}